\documentclass[12pt]{article}
\usepackage{amssymb}
\newtheorem{theorem}{Theorem}[section]
\newtheorem{proposition}[theorem]{Proposition}
\newtheorem{lemma}[theorem]{Lemma}
\newtheorem{definition}[theorem]{Definition}
\newtheorem{remark}[theorem]{Remark}
\newtheorem{corollary}[theorem]{Corollary}

\newcommand{\Irr}{\mathrm{Irr}}
\newcommand{\mod}{\mathrm{mod}}
\newcommand{\Sons}{\mathrm{Sons}}
\newcommand{\Father}{\mathrm{Father}}
\newcommand{\Des}{\mathrm{Des}}

\begin{document}
\title{Strongly regular Cayley graphs over primary abelian groups of rank 2}
\author{Yefim I. Leifman \\
Department of Mathematics and Computer Science, \\ Bar-Ilan University, \\ 52900, Ramat-Gan, Israel \\ email: leifmany@lycos.com
\\ \\
Mikhail E. Muzychuk \\
Department of Computer Science and Mathematics, \\ Netanya Academic College, \\ 42365, Netanya, Israel \\ email: muzy@netanya.ac.il}

\maketitle

{\bf Abstract. }
Strongly regular Cayley graphs with Paley parameters over abelian groups of rank 2 were studied in \cite{davis-94} and \cite{leung-95}. It was shown that such graphs exist iff the corresponding group is isomorphic to ${\mathbb Z}_{p^n} \oplus {\mathbb Z}_{p^n}$, where $p$ is an odd prime. In this paper we classify all strongly regular Cayley graphs over this group using Schur rings method. As a consequence we obtain a complete classification of strongly regular Cayley graphs with Paley parameters over abelian groups of rank 2.


\section{Introduction}

\begin{definition} [\cite{bose-63}] \label{d1.1}
An undirected graph without loops and multiple edges on $\nu$ vertices is called $(\nu,k,\lambda,\mu)$-strongly regular whenever there exist integers $k,\lambda,\mu$ satisfying 
\begin{enumerate}
\item each vertex is adjacent to $k$ other vertices,
\item each adjacent pair of vertices has $\lambda$ vertices, which are adjacent to both of them,
\item each non-adjacent pair of vertices has $\mu$ vertices, which are adjacent to both of them.
\end{enumerate}
\end{definition}

Let $R$ be a ring with unit 1. If $K$ is a subset of a finite group $G$, then the group-ring element $ \sum_{g \in K}{g}\in R[G]$ is called a {\it simple quantity} and will be denoted by $\underline K$.  An $R$-module with a basis $\{ \underline {T_1},\dots ,\underline {T_k} \} $ where $T_1,\dots , T_k$ are mutually disjoint sets with union $G$ is called an S-{\it module over} $G$ with a {\it standard basis} $\{ \underline {T_1},\dots ,\underline {T_k} \} $.

\begin{definition} [\cite{wielandt-64}] \label{d1.2}
An S-module $C$ over $G$ is called an S-ring over $G$ if the following conditions are satisfied
\begin{enumerate}
\item $C$ is a subring of $R[G]$,
\item $1\in C$,
\item if $\sum_{g\in G}{c_gg\in C}$, then $\sum_{g\in G}{c_gg^{-1}\in C}$.
\end{enumerate}
\end{definition}

$\Gamma_G(S)$ will denote a Cayley graph over a group $G$ with $S$ as a generating set. The following theorem is well known.
\begin{theorem}[\cite{ma-89}] \label{t1.3}
Let $S$ be a symmetric subset of a group $G$ (i.e. if $s  \in S$ then $s^{-1} \in S$) with $e \not\in S$. $\Gamma_G(S)$ is a strongly regular Cayley graph iff $\langle 1,\underline S,\underline G-\underline S-1 \rangle$ is an S-ring over $G$.
\end{theorem}

 Let $\Gamma_G(S)$ be a strongly regular Cayley graph (SRCG). If either $S\bigcup \{ e\} $ or $G\backslash S$ (the part of $G$ out of $S$) is a subgroup of $G$, then either $\Gamma_G(S)$ or its complement is a disjoint union of complete subgraphs of equal size. In this case we shall say that $\Gamma_G(S)$ is {\it trivial}. $\langle 1,\underline G \rangle $ is called a {\it trivial} S-ring over $G$.

\begin{definition} \label{d1.4}
An S-ring $C$ over $G$ is called primitive if $K=\{ e\}$ and $K=G$ are only subgroups of $G$ for which $\underline K\in C$ holds.
\end{definition}

It follows that the existence of a non-trivial SRCG over a given group $G$ implies the existence of a primitive S-ring over $G$. By Schur theorem \cite{wielandt-64} there is no non-trivial primitive S-ring over a cyclic group of composite order, by Kochendorfer's theorem \cite{kochendorfer-37} there is no non-trivial primitive S-ring over ${\mathbb Z}_{p^a} \oplus {\mathbb Z}_{p^b}$
with $a>b$ as well. A more general result \cite{leung-95} states that if a Sylow $p$-subgroup of the group is of type ${\mathbb Z}_{p^a} \oplus {\mathbb Z}_{p^b}$
with $a>b$, then there is no non-trivial primitive S-ring over this group. The structure of S-rings over ${\mathbb Z}_p$ is well known \cite{faradzev-90}. Each S-ring corresponds to a subgroup $H\leq Aut({\mathbb Z}_p)$: an S-ring has a standard basis of orbits of this subgroup. In particular, the 
SRCG or the Cayley tournament corresponds to the unique subgroup of $Aut({\mathbb Z}_p)$ of index 2 \cite{bridges-79}.

Therefore the "first family" of groups which is suitable for the search of non-trivial SRCGs is  ${\mathbb Z}_{p^n} \oplus {\mathbb Z}_{p^n}$ with $p$ prime. The SRCGs over ${\mathbb Z}_{p^n} \oplus {\mathbb Z}_{p^n}$ with {\it Paley parameters}, namely 
$(\nu ,(\nu -1)/2,(\nu -5)/4,(\nu -1)/4)$, were considered by J.A. Davis \cite{davis-94} $(n=2)$ and K.H. Leung, S.L. Ma \cite{leung-95} in the general case. Some examples of SRCGs were constructed in these papers but the enumeration problem was not considered. The goal of this paper is to describe all SRCGs over ${\mathbb Z}_{p^n} \oplus {\mathbb Z}_{p^n}$, where $p$ is an odd prime.

In order to formulate the main result, we need to introduce additional notations. Let $\Delta $ be the poset of all cyclic subgroups of ${\mathbb Z}_{p^n} \oplus {\mathbb Z}_{p^n}$ ordered by inclusion. The Hasse diagram of this poset is a tree with the trivial subgroup as a root. The valency of a node $H\in \Delta $ is 1 if $H$ is a leaf and $p+1$ otherwise. It turns out (Proposition \ref{p2.5}) that each SRCG over ${\mathbb Z}_{p^n} \oplus {\mathbb Z}_{p^n}$ is generated by a set of generators of elements of a subset of $\Delta $. We shall say that $S\subseteq \Delta $ defines an SRCG if $\Gamma _G(\bigcup _{H\in S}{O_H})$ is an SRCG, where $O_H$ is the set of all generators of $H$. We denote $\Gamma _G(\bigcup _{H\in S}{O_H})$ by $\Gamma _G(S)$  and the simple quantity $ \sum_{ H\in S}\underline {O_H}$  by $[S]$ for $S\subseteq \Delta $.

For a cyclic subgroup $H$ of ${\mathbb Z}_{p^n} \oplus {\mathbb Z}_{p^n}$ 
we define $\Father(H)=pH=\{ ph\mid h\in H\} $, $\Sons(H)=\{ F\in \Delta \mid \Father(F)=H\} $. If $|H| =p^i$, define the {\it length} of $H$ by $l(H)=i$.
Any set of the form $\Sons(H)$, $H\in \Delta $, will be called a {\it block} of $\Delta $. Any union of blocks will be called a {\it block set}. Two subsets $A_1,A_2\subseteq \Delta $ will be called {\it block equivalent} if their symmetric difference  $A_1\ominus A_2$ is a block set. 

\begin{definition}\label{d1.5}
Let $(a_1,\dots ,a_n)$ be an integer vector, $0\leq a_1\leq p+1$, $0\leq a_m\leq p-1$, $2\leq m\leq n$. We say that $S\subseteq \Delta $ is $(a_1,\dots ,a_n)$-homogeneous if $\{ e\} \not\in S $ and for each $H\in \Delta $ such that $0\leq l(H)<n$ it holds that 
\begin{displaymath}
|\Sons(H)\cap S|=
\left\{ 
\begin{array}{ll}
a_{l(H)+1}+1 & if \; H\in S, \\
a_{l(H)+1} & if \;  H\not\in S. 
\end{array}
\right.
\end{displaymath}
\end{definition}

A complement of a graph which is defined by an $(a_1,\dots ,a_n)$-homogeneous set is a graph which is defined by a $(p+1-a_1,p-1-a_2,\dots ,p-1-a_n)$-homogeneous set. We call two subsets $S,T\subseteq \Delta$ complement iff $S\cup T=\Delta \backslash \{\{e\}\}$ and $S\cap T=\emptyset $.

Let $\Delta _1$ be the set which contains all cyclic subgroups of ${\mathbb Z}_{p^n} \oplus {\mathbb Z}_{p^n}$ of order $p$ and the trivial subgroup. Let $S\subseteq \Delta $ define a non-trivial SRCG. There exists a unique homogeneous set $S^h\subseteq \Delta$ which is block equivalent to $S$ and satisfies $S^h\cap \Delta_1=S\cap \Delta_1$  (Corollary \ref{c5.9}).

The group $\langle (p^{n-1},0),(0,p^{n-1})\rangle $ is the group of all elements of order dividing $p$. The main theorem of this paper is:

\begin{theorem} \label{t1.6}
Let $p$ be a prime number. Every strongly regular Cayley graph over ${\mathbb Z}_{p^n} \oplus {\mathbb Z}_{p^n}$ is defined by a subset of $\Delta$. Let $p>2$ and let $\varphi :{\mathbb Z}_{p^n} \oplus {\mathbb Z}_{p^n}\rightarrow$ $ {\mathbb Z}_{p^n} \oplus {\mathbb Z}_{p^n}/ $ $\langle (p^{n-1},0),(0,p^{n-1})\rangle \cong $ ${\mathbb Z}_{p^{n-1}} \oplus {\mathbb Z}_{p^{n-1}}$ be the canonical homomorphism, $S\subseteq \Delta$ and $S\neq \emptyset$, $S\neq \Delta \backslash \{\{e\}\}$. $S$ defines a non-trivial strongly regular Cayley graph over ${\mathbb Z}_{p^n} \oplus {\mathbb Z}_{p^n}$ iff one of the following conditions is true:
\begin{enumerate}
\item $S$ is an $(a_1,a_2,\dots ,a_2)$-homogeneous set and $S$ is not a $(1,0,\dots ,0)$ or a $(p,p-1,\dots ,p-1)$-homogeneous set;
\item
if $n>3$, then $S^h$ is an $(a_1,0,\dots ,0,a_n)$-homogeneous set with $a_n>0$, $S^h\subseteq S$ and $Q=\varphi (S\backslash S^h)$ defines a non-trivial strongly regular Cayley graph over $\varphi (p({\mathbb Z}_{p^n} \oplus {\mathbb Z}_{p^n}))\cong {\mathbb Z}_{p^{n-2}} \oplus {\mathbb Z}_{p^{n-2}}$ for which $Q^h$ is a $(0,0,\dots ,0,a_n)$ or a $(p,p-1,\dots ,p-1,a_n-1)$-homogeneous set;

if $n=3$, then $S^h$ is an $(a_1,0,a_3)$-homogeneous set with $a_3>0$, $S^h\subseteq S$ and $Q=\varphi (S\backslash S^h)$ is an $(a_3)$-homogeneous set which defines a strongly regular Cayley graph over $\varphi (p({\mathbb Z}_{p^3} \oplus {\mathbb Z}_{p^3}))\cong {\mathbb Z}_{p} \oplus {\mathbb Z}_{p}$;
\item $S$ is a complement of the mentioned in the previous item.
\end{enumerate}
\end{theorem} 

All non-trivial SRCGs over ${\mathbb Z}_{p^n} \oplus {\mathbb Z}_{p^n}$ with $p>2$ are of the Latin Square Type with principal eigenvalue $k$ and non-principal eigenvalues $r$, $s$ such that $r=s+p^n$, $k=s-sp^n$. In the first case of the theorem $s=-a_1-a_2(p^n-p)/(p-1)$, where $0 \leq a_1 \leq p+1$, $0 \leq a_2 \leq p-1$ and $s \not\in \{ 0, -1, -p^n, -p^n-1 \}$. In two subcases of the second case 
$s=-a_1-a_np^{n-1}$, where $0 \leq a_1 \leq p+1$, $0 < a_n \leq p-1$. We do not consider the isomorphism problem of graphs with the same parameters in this paper.

\section{Strongly regular Cayley graphs and rational S-rings over finite abelian groups}

Let $G$ be a finite abelian group. Denote by $\Irr(G)$ the set of irreducible characters of $G$.
In what follows we extend an irreducible character of $G$ to the complex group-algebra ${\mathbb C}[G]$. But for simplicity of notations, if $\chi \in \Irr(G)$, then we write $\ker(\chi )$ instead of $\ker(\chi |_G)$, where $\chi |_G$ is a restriction of $\chi $ on $G$.
For $S\subseteq G$ and $t\in {\mathbb Z}$ we define $S^{(t)}=\{ g^t \mid g\in S\}$.
\begin{theorem} [\cite{leung-95}] \label{t2.1} 
Let $G$ be an abelian group of order $\nu$ and $S$ be a subset of $G$ with $e\not\in S$ and $S^{(-1)}=S$. Then $\Gamma_G(S)$ is an $(\nu ,k,\lambda ,\mu )$-SRCG over $G$ if and only if for any irreducible character $\chi $ of G
\begin{displaymath}
\chi (\underline S)=
\left\{ 
\begin{array}{ll}
k & if \; \chi \mbox{ is  principal on }G, \\
(\lambda -\mu\pm\sqrt{\delta })/2 & if \; \chi \mbox{ is  non-principal on } G, 
\end{array}
\right.
\end{displaymath}
where $\chi (\underline S)=\sum_{g\in S}{\chi (g)}$, $\delta =(\lambda-\mu)^2+4(k-\mu )$. These values are equal to the SRCG eigenvalues $k,r,s$ correspondingly.
\end{theorem}

\begin{theorem} [\cite{leung-95}] \label{t2.2}
Let $G$ be an abelian group of order $\nu $ and $S$ be a subset of $G$ with $e\not\in S$ and $S^{(-1)}=S$. Suppose that there exists an $(\nu ,k,\lambda ,\mu )$-SRCG $\Gamma_G(S)$ such that $\delta =(\lambda-\mu)^2+4(k-\mu )$ is not a square. Then $\Gamma_G(S)$ is an SRCG with Paley parameters $(\nu ,(\nu -1)/2,(\nu -5)/4,(\nu -1)/4)$ and $\nu =p^{2\eta +1}$ for some prime $p\equiv 1 \; \mod \; 4$ and integer $\eta $.
\end{theorem}

Let $G$ be a finite abelian group of exponent $m$. Let ${\bf P}(G)$ be the group consisting of all automorphisms of $G$ of the form $x\mapsto x^t$, where $t$ ranges through all residues $t$ which are relatively prime to $m$. The orbits of this action are in a one-to-one correspondence with cyclic subgroups of $G$. More precisely, if $H\leq G$ is a cyclic subgroup, then the set of its generators $O_H$ is an orbit of ${\bf P}(G)$. Denote the S-module (with $\mathbb C$ as $R$) with standard basis of simple quantities $\underline{O_H}$ by $W(G)$. 

\begin{theorem} [\cite{bridges-82}] \label{t2.3} 
Let $G$ be a finite abelian group. Then the S-module $W(G)$ is an S-ring over $G$. Moreover, $W(G)$ is the unique maximal S-ring over $G$ for which the values of the irreducible characters of $G$ on the elements of its standard basis are rational.
\end{theorem}

\begin{definition}\label{d2.4}
Let $G$ be a finite abelian group. Any S-ring over $G$ contained in $W(G)$ is called a rational S-ring over $G$.
\end{definition}

\begin{proposition}\label{p2.5}
Let $A$ be a finite abelian group of order which is not of the form $p^{2\eta +1}$.
There exists a one-to-one correspondence between the following sets:
\begin{enumerate}
\item rank 3 rational S-rings over $A$,
\item pairs of complement SRCGs over $A$,
\item pairs of complement unions of standard basic subsets of $W(A)$, excepting $\{ e\} $, for which the set of all values of non-principal irreducible characters of $A$ on each union contains only two elements.
\end{enumerate}
\end{proposition}

\begin{definition} [\cite{bridges-82}] \label{d2.6}
Let $\lambda :G\times G\rightarrow {\mathbb C}^*$ satisfy 
\begin{enumerate}
\item $\lambda (g,h)=\lambda (h,g)$,
\item $\lambda (g,h_1h_2)=\lambda (g,h_1)\lambda (g,h_2)$,
\item $\forall g\in G\; \lambda (g,h)=1 \; \Rightarrow \; h=e$.
\end{enumerate}
Then $g\mapsto \lambda (g,-)$ is called a symmetric isomorphism of $G$ with its character group.
\end{definition}

\begin{definition} \label{d2.7}
Let $\Gamma_G(S)$ be an SRCG and $\lambda $ be a symmetric isomorphism of $G$ with its character group. Define $S^+$ such that $e\not\in S^+$ and for each $g\in G$, $g\neq e$ it holds that $g\in S^+$ iff $\sum_{h \in S}{\lambda (g,h)}=r$, where $r$ is the largest non-principal eigenvalue of  $\Gamma_G(S)$. Then $\Gamma_G(S^+)$ is called the dual graph to $\Gamma_G(S)$ with respect to $\lambda $.
\end{definition}

\begin{theorem} [\cite{delsarte-73}] \label{t2.8}
$\Gamma_G(S^+)$ is a non-trivial SRCG iff $\Gamma_G(S)$ is a non-trivial SRCG and in this case $(r-s)(r^+-s^+)=|G|$, where $r^+$, $s^+$ are non-principal eigenvalues of $\Gamma_G(S^+)$.
\end{theorem}

\section{Complex characters of rational S-rings over \\ finite abelian groups}

\begin{proposition} [\cite{bridges-79}] \label{p3.1}
The set of simple quantities which correspond to cyclic subgroups of $G$ forms a basis of $W(G)$ called the subgroup basis. 
\end{proposition}
{\bf Proof. } Let $C_m$ be a cyclic subgroup of $G$ of order $m$, $\underline {C_m}=\sum_{l\in D_m}\underline {O_l}$, where $D_m$ is the set of all divisors of $m$, $O_l$ is the orbit corresponding to the cyclic subgroup $C_l$ of $C_m$. Then $\underline {O_m}=\sum_{l\in D_m}{\mu (m/l)\underline {C_l}}$, where $\mu (x)$ is the M{\"o}bius function. $\square $

\begin{proposition} \label{p3.2}
Let $\rho ,\sigma \in \Irr(G)$. $\rho $ and $\sigma $ are equal on $W(G)$ iff 
$\ker(\rho ) =\ker(\sigma )$.
\end{proposition}
{\bf Proof. } Let $\underline H$ be an element of the subgroup basis of $W(G)$.
\begin{displaymath}
\sigma (\underline H)=\sum_{h\in H}{\sigma(h)}=
\left\{ 
\begin{array}{ll}
0 & \mbox{if} \; H\not\subseteq \ker(\sigma ), \\
|H| & \mbox{if} \; H\subseteq \ker(\sigma ). 
\end{array}
\right.
\end{displaymath}

Thus if $\ker(\sigma )=\ker(\rho )$, then $\rho$ and $\sigma $ are equal on the subgroup basis of $W(G)$. If $\ker(\sigma ) \neq \ker(\rho )$, then there exists $h\in G$ for which $h\in \ker(\sigma )$ and $h\not\in \ker(\rho )$ or $h\not\in \ker(\sigma )$ and $h\in \ker(\rho )$. Assume that $h\in \ker(\sigma )$ and $h\not\in \ker(\rho )$. Then $H=\langle h\rangle \subseteq \ker(\sigma )$ and  $H\not\subseteq \ker(\rho )$. Then $\rho (\underline H)\neq \sigma (\underline H)$ as desired. $\square$

\begin{definition} \label{d3.3}
A subgroup $H\leq G$ is called a cocyclic subgroup iff $G/H$ is a cyclic group.
\end{definition}

\begin{proposition} \label{p3.4}
There exists a one-to-one correspondence between the set of equivalence classes of $\Irr(G)$, where two characters belong to the same class iff they are equal on $W(G)$, and the set of cocyclic subgroups of $G$.
\end{proposition}

\begin{definition} \label{d3.5}
The intersection of all maximal subgroups of a group $G$ is called the Frattini subgroup of $G$ and denoted by $\Phi (G)$. If $G$ has no maximal subgroups, $\Phi (G)=G$ by definition.
\end{definition}

If $G$ is a cyclic group of the order $p_1^{a_1}\cdots p_s^{a_s}$, then $\Phi (G)$ has the index $p_1\cdots p_s$.

\begin{lemma} \label{l3.6}
Let $\chi \in \Irr(G)$, $h\in G$, $H=\langle h\rangle$, $F=H\cap \ker(\chi )$. Denote by $O_H$ the set of all generators of $H$. Then $$\chi (\underline {O_H})=|\Phi (H)|\mu \left( {|H|} \over {|F|}\right)
\varphi \left( {|F|} \over {|\Phi (H)|}\right).$$
\end{lemma}
{\bf Proof. } Denote by $o(h)$ the order of $h$. Since $H=\bigcup_{d|o(h)}{O_{\langle h^d\rangle}}$, we have $\chi (\underline H)=\sum_{d|o(h)}{\chi (O_{\langle h^d\rangle })}$. Using the M{\"o}bius inversion we obtain
 $$\chi (\underline {O_H})=\sum_{d|o(h)}{\mu (o(h)/d)\chi (\underline {G_d})},$$ where $G_d$ is the unique subgroup of $H$ of order $d$. By Proposition \ref{p3.2} we obtain $\chi (\underline {G_d})=0$ if  $G_d\not\subseteq \ker(\chi )$ and $\chi (\underline {G_d})=|G_d|$ if $G_d\subseteq \ker(\chi )$.
Thus $$\chi (\underline{O_H})=\sum_{d\bigl| |F|}{\mu (o(h)/d)\chi (\underline {G_d})}=\sum_{d\bigl| |F|}{\mu (o(h)/d)|G_d|}=\sum_{d\bigl| |F|}{\mu (o(h)/d)d}.$$
Since $o(h)/d=|H|/d=(|H|/|F|)(|F|/d)$ and $\mu (x)=0$ whenever $x$ is divisible by a square of a prime integer, $\chi (\underline {O_H})=0$ in the case of $|H|/|F|$ being divisible by a square of a prime integer.

If $|H|/|F|$ is not divisible by a square, then $F\supseteq \Phi (H)$ and
$$\chi (\underline {O_H})=|\Phi (H)|\sum_{|\Phi (H)|\bigl| d\bigr| |F|}{\mu \left( {o(h)} \over {d}\right)
\left( {d} \over {|\Phi (H)|}\right)}.$$
We have $o(h)/d=(o(h)/|F|)(|F|/d)$. Since $o(h)/d$ is not divisible by a square, $o(h)/|F|$ and $|F|/d$ are co-prime. Taking into account that $\mu $ is multiplicative in this case we obtain
$$\chi (\underline {O_H})=|\Phi (H)|\mu \left( {o(h)} \over {|F|}\right)\sum_{|\Phi (H)|\bigl| d\bigr| |F|}{\mu \left( {|F|} \over {d}\right)
\left( {d} \over {|\Phi (H)|}\right)}=$$
$$=|\Phi (H)|\mu \left( {|H|} \over {|F|}\right)\sum_{d_1\bigl| (|F|/|\Phi (H)|)}{\mu \left( {|F|} \over {|\Phi (H)|d_1}\right)d_1}=$$
$$=|\Phi (H)|\mu \left( {|H|} \over {|F|}\right)
\varphi \left( {|F|} \over {|\Phi (H)|}\right).$$\hfill $\square$

\section{The characters of the S-ring $W({\mathbb Z}_{p^n} \oplus {\mathbb Z}_{p^n})$}

In the following sections $G$ will stand for ${\mathbb Z}_{p^n} \oplus {\mathbb Z}_{p^n}$.

\begin{proposition} \label{p4.1}
The subgroups
$$\{ H\mid H=\langle (p^m,ap^m)\rangle, 0\leq m\leq n-1,0\leq a\leq p^{n-m}-1\} $$
$$\cup \{ H\mid H=\langle (bp^{m+1},p^m)\rangle ,0\leq m\leq n-1, 0\leq b\leq p^{n-m-1}-1\} \cup \{ \{ (0,0)\} \}$$
exhaust the set of cyclic subgroups of $G$. The set of cyclic subgroups is partially ordered by inclusion. The Hasse diagram of this poset is a tree, where the trivial subgroup is the root.
\end{proposition}
{\bf Proof. } The above subgroups are distinct and a total number of their generators is equal to $|G|$. Therefore they exhaust the set of cyclic subgroups of $G$. The lattice of subgroups of a cyclic $p$-group is a chain, therefore the Hasse diagram of the poset of cyclic subgroups of a $p$-group is a tree. This tree is denoted by $\Delta $.  $\square$

\begin{proposition} \label{p4.2} $\;$
\begin{enumerate}
\item
The subgroups
$$\{ K \mid K=\langle (1,a),(0,p^m)\rangle, 1\leq m\leq n,0\leq a\leq p^m-1\} $$
$$\cup \{ K \mid K=\langle (bp,1),(p^m,0)\rangle ,1\leq m\leq n, 0\leq b\leq p^{m-1}-1\} \cup \{ G \}$$ exhaust the set of cocyclic subgroups of $G$.
\item
Let $H\cong {\mathbb Z}_{p^n} \oplus {\mathbb Z}_{p^m}$ be a cocyclic subgroup of $G$.
Define $H^\Delta =\{ p^mh \mid h\in H\} $. Then $H\mapsto H^\Delta$ is a bijection between the set of cocyclic subgroups of $G$ and the set of cyclic subgroups of $G$. Moreover, $H_1\subseteq H_2$ iff $H_1^\Delta \supseteq H_2^\Delta$.
\item
The set of cocyclic subgroups is partially ordered by inclusion. The Hasse diagram of this poset is a tree, where $G$ is the root.
\end{enumerate}
\end{proposition}

{\bf Proof. } (1) If $H$ is a cocyclic subgroup of $G$, then $H\cong {\mathbb Z}_{p^n} \oplus {\mathbb Z}_{p^m}$ with $0\leq m\leq n$. Then the conclusion follows from Proposition \ref{p4.1}.

(2) $\langle (1,a),(0,p^{n-m})\rangle ^\Delta =\langle (p^m,ap^m)\rangle$ are distinct for $0\leq a\leq p^{n-m}-1$, $0\leq m\leq n$.  $\langle (1,a_1),(0,p^{m_1})\rangle \supseteq \langle (1,a_2),(0,p^{m_2})\rangle \Leftrightarrow (a_1 \equiv a_2(\mod \; p^{m_1})$ and $m_1\leq m_2) \Leftrightarrow \langle (p^{n-m_1},a_1p^{n-m_1})\rangle \subseteq \langle (p^{n-m_2},a_2p^{n-m_2})\rangle$. $\square$

Denote the tree of cocyclic subgroups by $\nabla $.
Thus $^\Delta $ maps the Hasse diagram of the poset $\nabla$ onto the Hasse diagram of $\Delta$. We denote the inverse function of $^\Delta$ by $^\nabla$.

It is easy to see that $|H|=p^n$ iff $H^\Delta =H$. This fact is generalized in the following lemma.

\begin{lemma} \label{l4.3}
Let $F\in \Delta$ and $H\in \nabla$. Then $F\subseteq H$ iff $l(F)-l(F\cap H^\Delta)\leq n-l(H^\Delta)$.
\end{lemma}
{\bf Proof. } Let $|H|=p^{n+m}$. Then $|H^\Delta |=p^{n-m}$, $l(H^\Delta)=n-m$ . So our claim is equivalent to $F\subseteq H$ iff $l(F)-l(F\cap H^\Delta)\leq m$. Assume $F\subseteq H$. Then $H^\Delta =p^mH\supseteq p^mF$. Therefore $F\cap H^\Delta \supseteq p^mF$. $F$ is cyclic and $[F:p^mF]\leq p^m$ hence $[F:F\cap H^\Delta ]\leq [F:p^mF]\leq p^m$. But $[F:F\cap H^\Delta ]=|F|/|F\cap H^\Delta |=p^{l(F)-l(F\cap H^\Delta )}$ implies $l(F)-l(F\cap H^\Delta )\leq m$, as desired.

If $l(F)<m$, then $F\subseteq H$. Assume now that $l(F)\geq m$ and $l(F)-l(F\cap H^\Delta )\leq m$, i.e $[F:F\cap H^\Delta ]\leq p^m=[F:p^mF]$. Since $F$ is cyclic, this inequality is equivalent to the inclusion $F\cap H^\Delta \supseteq p^mF$ which is equivalent to $H^\Delta \supseteq p^mF$. We may assume that $H=\langle (1,0),(0,p^{n-m})\rangle$. Then $H^\Delta=\langle (p^m,0)\rangle$. Let $f=(f_1,f_2)$ be a generator of $F$. Then $p^mF=\langle (p^mf_1,p^mf_2)\rangle$. Now $p^mF\subseteq H^\Delta$ implies $f_2\equiv 0(\mod \; p^{n-m})$. Therefore $f=(f_1,p^{n-m}f'_2)\in 
\langle (1,0),(0,p^{n-m})\rangle =H$. $\square$

Now we extend some notations from the previous sections. Define $\Delta_i=\{H\in \Delta \mid l(H)\leq i\}$. Define {\it descendants} of $H\in \Delta $ inductively: $\Des_1(H)=\Sons(H)$, $\Des_{i+1}(H)=\{ F\in \Delta \mid \Father(F)\in \Des_i(H)\}$. Let $\Delta _{i}^{j}=\{ \Des_j(H) \mid H\in \Delta_i\}$ for $i+j\leq n$. We write $\Delta ^j$ instead of $\Delta _{n-j}^j$. $\Delta _{i}^{j}$ is a poset with order relation $"\subseteq ^j"$: $F_1\subseteq^jF_2$, $F_1=\Des_j(H_1)$, $F_2=\Des_j(H_2)\in \Delta _{i}^{j}$ iff $H_1\subseteq H_2$. The Hasse diagram of this poset is a tree. Similarly to $\Delta $, we can define for $\Delta _{i}^{j}$ blocks, the block equivalence and the functions  $l(F)$, $\Sons(F)$, $\Des_m(F)$ for $F\in \Delta _{i}^{j}$. For $T\subseteq \Delta$ we define $T^j=\{ F\in \Delta^j \mid F\subseteq T\}$. For $T\subseteq \Delta _{i}^{j}$ we define $T_m=T\cap \Delta _{m}^{j}$. Denote $T_i^j=(T^j)_i$.

Similarly to $\Delta $, for $\nabla $ we can define $\nabla_i$, a block in $\nabla_i$ and $l(F)$, $\Sons(F)$, $\Des_m(F)$ for $F\in \nabla_i$.

If $S\subseteq \Delta _{i}^{j}$, then $[S]$ will mean the simple quantity 
$\sum_{F\in \bigcup _{H\in S}H}\underline {O_F}$. In what follows the notation $\chi _H$ will mean an irreducible character with a kernel $H$ and we shall write $\chi _H[S]$ instead of $\chi _H([S])$.

\begin{corollary} \label{c4.4}
Let $F\in \Delta$ and $H\in \nabla$. Then
\begin{displaymath}
\chi_H[F]=\chi_H(\underline {O_F})=
\left\{ 
\begin{array}{ll}
|F|(p-1)/p & if \;\; l(F)-l(F\cap H^\Delta)\leq n-l(H^\Delta), \\
-|F|/p & if \;\; l(F)-l(F\cap H^\Delta)= n-l(H^\Delta)+1, \\ 
0 & otherwise.
\end{array}
\right.
\end{displaymath}
\end{corollary}
{\bf Proof. } Apply Lemma \ref{l4.3} and Lemma \ref{l3.6}. $\square$

\begin{proposition} \label{p4.5}
Let $X,Y\subseteq \Delta_i^j$ such that $(X\cup Y)\cap \Delta_0^j=\emptyset$ and for each $l$,  $1\leq l\leq i$, either $X\cap (\Delta_l^j \backslash \Delta_{l-1}^j)=\emptyset$ or $Y\cap (\Delta_l^j \backslash \Delta_{l-1}^j)=\emptyset$. Then for each $H_1,H_2\in \nabla_{n-j}\backslash \nabla_{n-1-j}$ it holds that
$$|(\chi_{H_1}[X]-\chi_{H_1}[Y])-(\chi_{H_2}[X]-\chi_{H_2}[Y])|\leq (2p^i-1)p^{2j}.$$
\end{proposition}
{\bf Proof. } If $H\in \nabla_{n-j}\backslash \nabla_{n-1-j}$, then
$-p^{2j+1}\leq \chi_H[X\cap (\Delta_1^j\backslash \Delta_0^j)]\leq (p-1)p^{2j}$ and 
$-(p^m-p^{m-1})p^{2j}\leq \chi_H[X\cap (\Delta_m^j\backslash \Delta_{m-1}^j)]\leq (p^m-p^{m-1})p^{2j}$, $2\leq m\leq i$. Summarizing these inequalities we have $-p^{i+2j}\leq \chi_H[X]-\chi_H[Y]\leq (p^i-1)p^{2j}$ if $Y\cap (\Delta_1^j \backslash \Delta_0^j)=\emptyset$ and $-(p^i-1)p^{2j}\leq \chi_H[X]-\chi_H[Y]\leq p^{i+2j}$ if $X\cap (\Delta_1^j \backslash \Delta_0^j)=\emptyset$. Writing these inequalities for $H_1$ and for $H_2$ we have the claim. $\square$

\section{Homogeneous strongly regular Cayley graphs over ${\mathbb Z}_{p^n} \oplus {\mathbb Z}_{p^n}$}

\begin{proposition} [\cite{brouwer-89}] \label{p5.1} 
Let $\Gamma$ be a strongly regular graph. If one of its  eigenvalues is $0$ or $-1$, then $\Gamma$ is a trivial strongly regular graph.
\end{proposition}

\begin{proposition} \label{p5.2}
An SRCG over $G$ defined by a block subset $S\subseteq \Delta$ is trivial.
\end{proposition}
{\bf Proof. } Let $H\in \nabla\backslash \nabla_{n-1}$. Corollary \ref{c4.4} and $H=H^\Delta$ imply, that if $B$ is a block $\Delta_1\backslash \Delta_0$, then $\chi_H[B]=-1$ and if a block $B\subseteq \Delta_i\backslash \Delta_{i-1}$, $2\leq i\leq n$, then $\chi_H[B]=0$. Therefore $\chi_H[S]\in \{ 0,-1\}$, whence by Proposition \ref{p5.1} $S$ defines a trivial SRCG. 
 $\square$

\begin{proposition} \label{p5.3}
If $S\subseteq \Delta$ is not a block set, then there exist $m$, $1\leq m\leq n$, and $H_1,H_2\in \nabla\backslash \nabla_{n-1}$ such that $\chi_{H_1}[S]-\chi_{H_2}[S]=p^m$.
\end{proposition}
{\bf Proof. } Let $m$ be a maximal number for which there exists a block $B\subseteq \Delta_m\backslash \Delta_{m-1}$ such that $B\cap S\neq \emptyset$, $B\backslash S\neq \emptyset $. Let $F_1\in B\cap S$, $F_2\in B\backslash S$. We set $H_1$ and $H_2$ such that $H_1^\Delta \in \Des_{n-m}(F_1)$ and $H_2^\Delta \in \Des_{n-m}(F_2)$. $\square$

\begin{proposition} \label{p5.4}
If $\Gamma_G(S)$ is a non-trivial SRCG with non-principal eigenvalues $r$ and $s$, then $r-s=p^n$. In particular, $\Gamma_G(S)$ is of Latin or Negative Latin Square Type.
\end{proposition}
{\bf Proof. } By Propositions \ref{p5.2} and \ref{p5.3} $r-s=p^m$, $1\leq m\leq n$. By Theorem \ref{t2.8} an SRCG $\Gamma_G(S)$ is non-trivial iff a dual SRCG $\Gamma_G(S^+)$ is non-trivial. Therefore $(r-s)(r^+-s^+)=|G|=p^{2n}$ implying $r-s=p^n$.

Each $(\nu ,k,\lambda ,\mu )$-strongly regular graph satisfies the equality $(\nu -k-1)\mu=k(k-1-\lambda)$ \cite{seidel-69}. Therefore by Theorem \ref{t2.1} $(k-r)(k-s)/\nu=\mu$ and $k+rs=\mu$, from which it follows that $k+rs=(k-r)(k-s)/\nu$. Moreover, every non-trivial SRCG over $G$ satisfies $r-s=p^n$. Therefore either its valency is equal to $k'=s-sp^n$ or $k''=s+sp^n+p^n+p^{2n}$. A $(\nu ,k,\lambda ,\mu )$-strongly regular graph is called a Latin Square Type strongly regular graph if its parameters are of the form $\nu =m^2$, $k=q(m-1)$, $\lambda =m-2+(q-1)(q-2)$, $\mu =q(q-1)$. It is called a Negative Latin Square Type strongly regular graph if its parameters are of the form $\nu =m^2$, $k=q(m+1)$, $\lambda =-m-2+(q+1)(q+2)$, $\mu =q(q+1)$. In the case of ${\mathbb Z}_{p^n} \oplus {\mathbb Z}_{p^n}$ $k=k'$ implies that the graph is of Latin Square Type and $k=k''$ implies that the graph is of Negative Latin Square Type. $\square$

\begin{definition} \label{d5.5}
Let $(a_1,\dots ,a_i)$ be an integer vector, $0\leq a_1\leq p+1$, $0\leq a_m\leq p-1$, $2\leq m\leq i$. We say that $S\subseteq \Delta_i^j $ is $(a_1,\dots ,a_i)$-homogeneous if $\Delta_0^j \not\subseteq S $ and for each $H\in \Delta_{i-1}^j $ it holds that 
\begin{displaymath}
|\Sons(H)\cap S|=
\left\{ 
\begin{array}{ll}
a_{l(H)+1}+1 & if \; H\in S, \\
a_{l(H)+1} & if \; H\not\in S. 
\end{array}
\right.
\end{displaymath}
\end{definition}

Let $S\subseteq \Delta_k^j $ be an  $(a_1,\dots ,a_k)$-homogeneous set. Fix $F\in \Delta_k^j$, $F\not\in S$. If $l(F)=0$, then $|\Des_1(F)\cap S|=a_1$ and 
$$|\Des_{m+1}(F)\cap S|=(a_{m+1}+1)|\Des_m(F)\cap S|+a_{m+1}(p^m+p^{m-1}-|\Des_m(F)\cap S|).$$ By induction $|\Des_m(F)\cap  S|=a_1+\sum_{i=2}^{m}{a_i(p^{i-1}+p^{i-2})}$. Analogously,  if $l(F)\neq 0$, then $|\Des_{m-l(F)}(F)\cap S|=\sum_{i=l(F)+1}^{m}{a_ip^{i-l(F)-1}}$. Since these numbers depend only on $l(F)$, we set $A_{l(F),m}=|\Des_{m-l(F)}(F)\cap S|$. As it was shown before $A_{0,m}=a_1+\sum_{i=2}^{m}{a_i(p^{i-1}+p^{i-2})}$, $A_{l,m}=\sum_{i=l+1}^{m}{a_ip^{i-l-1}}$. Analogously, if $F\in \Delta_k^j$, $F\in S$, then $|\Des_{m-l(F)}(F)\cap S|=A_{l(F),m}+1$. Finally, $|\Des_{m-l(F)}(F)\cap S|=A_{l(F),m}+\delta_S(F)$, where $\delta_S(F)=1$ if $F\in S$ and $\delta_S(F)=0$ if $F\not\in S$. In particular, $A_{m,m+1}=a_{m+1}$. Also we set $A_{m,m}=0$.

\begin{proposition} \label{p5.6}
Let $S\subseteq \Delta$ be an $(a_1,\dots ,a_n)$-homogeneous set. Then for all $t$, $1\leq t \leq n $
\begin{enumerate}
\item 
$|\{\chi _H[S]\mid l(H)=t\}|\leq 2$,
\item
$\chi_{H_1}[S]\equiv \chi_{H_2}[S](\mod \; p^n)$ for all $H_1,H_2\in \nabla_t\backslash \nabla_{t-1}$,
\item 
if $(a_1,\dots ,a_t)\neq (0,\dots ,0)$ and $(a_1,\dots ,a_t)\neq (p+1,p-1,\dots ,p-1)$, then there exist $X,Y\in \nabla_t\backslash \nabla_{t-1}$ such that $\chi_X[S]-\chi_Y[S]=p^n$,
\item
if $(a_1,\dots ,a_t)= (0,\dots ,0)$ or  $(a_1,\dots ,a_t)=(p+1,p-1,\dots ,p-1)$, then for all $X,Y\in \nabla_t\backslash \nabla_{t-1}$ it holds that $\chi_X[S]=\chi_Y[S]$.
\end{enumerate}
\end{proposition}
{\bf Proof. } By Corollary \ref{c4.4} each character partitions $\Delta $ into three subsets. Let $H\in \nabla_t\backslash \nabla_{t-1}$ and $S=S'_H\cup S''_H\cup S'''_H$ be a partition defined as follows: $$S'_H=\{ F\in S\mid \chi_H[F]=|F|(p-1)/p\},$$
$$ S''_H=\{ F\in S\mid \chi_H[F]=-|F|/p\},$$
$$ S'''_H=\{ F\in S\mid \chi_H[F]=0\}.$$
Then $\chi_H[S]= \chi_H[S'_H]+ \chi_H[S''_H]$. Define $F_t=H^\Delta$, $F_{i-1}=\Father(F_i)$, $i=t,\dots ,2$. Then
$$\chi _H[S'_H]=\sum_{i=1}^{n-t}{A_{0,i}(p^i-p^{i-1})}+
\sum_{i=1}^{t}{(A_{i,n-t+i}+\delta_S(F_i))(p^{n-t+i}-p^{n-t+i-1})},$$
$$\chi_H[S''_H]=-\sum_{i=0}^{t-1}{(A_{i,n-t+i+1}+\delta_S(F_i)-A_{i+1,n-t+i+1}-\delta_S(F_{i+1}))p^{n-t+i}}.$$
Adding these equalities we obtain 
$$\chi _H[S]=\sum_{i=1}^{n-t}{A_{0,i}(p^i-p^{i-1})}-A_{0,n-t+1}p^{n-t}$$
\begin{equation} \label{e5.6.1}
+\sum_{i=1}^{t-1}{(A_{i,n-t+i}-A_{i,n-t+i+1})p^{n-t+i}}+(A_{t,n}+\delta_S(H^\Delta))p^n.
\end{equation}
If  $(a_1,\dots ,a_t)\neq (0,\dots ,0)$ and $(a_1,\dots ,a_t)\neq (p+1,p-1,\dots ,p-1)$, then there exists a block $B\subseteq \nabla_t\backslash \nabla_{t-1}$ such that $B^\Delta \cap S\neq \emptyset$, $B^\Delta \backslash S\neq \emptyset$. We set $X\in B\cap S^\nabla$, $Y\in B \backslash S^\nabla$. $\square$

Let $S\subseteq \Delta$ or $S\subseteq \Delta_i^j$. Denote $x_m[S]=\min\{ \chi _H[S]\mid H\in \nabla_m
\backslash \nabla_{m-1}\}$ (all considered characters have rational values).

\begin{proposition} \label{p5.7}
Let $S\subseteq \Delta_i^j$, $\Delta_0^j\cap S=\emptyset$. If  $\chi_H[S]\equiv \chi_{H'}[S](\mod \; p^{i+2j})$ for every $H,H'\in \nabla_{n-j}\backslash \nabla_{n-j-1}$, then
\begin{enumerate}
\item $|\{ \chi_H[S]\mid H\in \nabla_{n-j}\backslash \nabla_{n-j-1}\} |\leq 2$, 
\item there exists a unique homogeneous set $S^h\subseteq \Delta_i^j$ which is block equivalent to $S$ and satisfies $x_{n-j}[S]=x_{n-j}[S^h]$;
it holds that $\chi _H[S]=\chi _H[S^h]$ for all $H\in \nabla_{n-j}\backslash \nabla_{n-j-1}$,
\item $x_{n-j}[S]=-p^{2j}\sum_{l=1}^{i}{a_lp^{l-1}}$ whenever $S$ is an $(a_1,a_2,\dots ,a_i)$-homogeneous set and $(a_1,a_2,\dots ,a_i)\neq (p+1,p-1,\dots ,p-1)$.
\end{enumerate}
\end{proposition}
{\bf Proof. } (1) By Proposition \ref{p4.5} $|\chi_H[S]- \chi_{H'}[S]|\leq (2p^i-1)p^{2j}$ whenever $H,H'\in \nabla_{n-j}\backslash \nabla_{n-j-1}$.

(2) We shall prove the claim by induction on $i$ with fixed $j$.

In the case of $i=1$ we set $S^h=S$ and this is the unique possibility to satisfy the condition $x_{n-j}[S^h]=x_{n-j}[S]$ since $x_{n-j}[\Delta_1^j\backslash \Delta_0^j]\not\equiv 0(\mod \; p^{2j})$.

Assume now that $i>1$. Let $H\in \nabla_{n-j}\backslash \nabla_{n-j-1}$. Consider the set $ S_{i-1}=S\cap \Delta_{i-1}^j$. By Corollary \ref{c4.4}, $\chi_H[F]\equiv 0(\mod \; p^{i-1+2j})$ whenever $l(F)=i$, therefore  $\chi_H[S_{i-1}]\equiv \chi_{H}[S](\mod \; p^{i-1+2j})$. Then by induction hypothesis $S_{i-1}$ is block equivalent to the unique  homogeneous set $(S_{i-1})^h$ which satisfies $x_{n-j}[S_{i-1}]=x_{n-j}[(S_{i-1})^h]$. Since $\chi_H[B]=-p^{2j}$ for a block $B=\Delta_1^j\backslash \Delta_0^j$ and $\chi_H[B]=0$ for a block $B\subseteq \Delta_k^j\backslash \Delta_{k-1}^j$, $k>1$, it holds that $(S_{i-1})^h\cap \Delta_1^j= S_{i-1}\cap \Delta_1^j$ and $\chi_H[S_{i-1}]=\chi_H[(S_{i-1})^h]$.

Denote $F=\Des_j(F^*)$, where $F^*$ is the unique forefather  of $H^\Delta$ of length $i$. Set $F'=\Father(F)$, $B=\Sons(F')$ (i.e. $B$ is the unique block that contains $F$). Then 
$$\chi_H[S]= \chi_H[S_{i-1}]+\chi_H[S\backslash \Delta_{i-1}^j]=\chi_H[(S_{i-1})^h]+(-|S\cap B|p^{i-1+2j}+p^{i+2j}\delta_S(F))$$
$$=x_{n-j}[(S_{i-1})^h]+\delta_{(S_{i-1})^h}(F')p^{i-1+2j}-|S\cap B|p^{i-1+2j}+p^{i+2j}\delta_S(F).$$
Since all $ \chi_H[S]$, $ H\in \nabla_{n-j}\backslash \nabla_{n-j-1}$ have the same residue modulo $p^{i+2j}$ by assumption, there exists $a_i\in [0,p-1]$ such that $|S\cap B|-\delta_{(S_{i-1})^h}(F')\equiv a_i(\mod \; p)$. The left-hand side belongs to  $[-1,p]$. If $0<a_i<p-1$, then $a_i=|S\cap B|-\delta_{(S_{i-1})^h}(F')$ and $(S\backslash \Delta_{i-1}^j)\cup (S_{i-1})^h$ is the homogeneous set we are looking for. If $a_i=0$, then there are three possibilities: $|S\cap B|=0$ and $\delta_{(S_{i-1})^h}(F')=0$ or $|S\cap B|=p$ and  $\delta_{(S_{i-1})^h}(F')=0$ or $|S\cap B|=1$ and $\delta_{(S_{i-1})^h}(F')=1$. Thus we can obtain $S^h$ by removing all blocks with $|S\cap B|=p$, $\delta_{(S_{i-1})^h}(F')=0$. Analogously, if $a_i=p-1$, we can obtain $S^h$ by adding all blocks with $|S\cap B|=0$, $\delta_{(S_{i-1})^h}(F')=1$.

$(S^h)_{i-1}$ is a homogeneous set and meets the assumption, therefore $(S^h)_{i-1}=(S_{i-1})^h$. Then $S^h$ is the unique homogeneous set for which $\chi_H[S^h]=\chi_H[S]$ for all $H\in \nabla_{n-j}\backslash \nabla_{n-j-1}$ by construction.

(3) The equation is a straightforward consequence of Corollary \ref{c4.4}. $\square$

\begin{corollary} \label{c5.8}
If $S\subseteq \Delta_i^j$, $S \cap \Delta_0^j=\emptyset$ and $S^h$ is an $(a_1,\dots ,a_i)$-homogeneous set, then
\begin{enumerate}
\item if $0<a_l<p-1$, then $S\cap (\Delta_l^j\backslash \Delta_{l-1}^j)=S^h\cap (\Delta_l^j\backslash \Delta_{l-1}^j)$ for each $l$, $1\leq l\leq i$,
\item if $S^h$ is neither $(0,\dots ,0)$-homogeneous nor $(p+1,p-1,\dots ,p-1)$-homogeneous, then there exists a block $B\subseteq \Delta_i^j\backslash \Delta_{i-1}^j$ such that $B\cap S\neq \emptyset$ and $B\backslash S\neq \emptyset$.
\end{enumerate}
\end{corollary}

\begin{corollary} \label{c5.9}
If $S\subseteq \Delta$ defines a non-trivial SRCG over $G$, then there exists a unique $(a_1,\dots ,a_n)$-homogeneous set $S^h\subseteq \Delta$ such that $S$ and $S^h$ are block equivalent and $s=x_n[S^h]$.
\end{corollary}

\begin{proposition} \label{p5.10}
Let $S\subseteq \Delta$ be an $(a_1,\dots ,a_n)$-homogeneous set. Then $S$ defines an SRCG iff $a_2=\cdots =a_n$.
\end{proposition}
{\bf Proof. }  
$$\chi_G[S_m]=\sum_{i=1}^{m}{A_{0,i}(p^i-p^{i-1})}=a_1(p^m-1)+\sum_{i=2}^{m}{a_i(p^{i-1}+p^{i-2})(p^m-p^{i-1})}.$$
If $H\in \nabla_t\backslash \nabla_{t-1}$, $1\leq t<n$, then $\chi_H[S_{n-t}]=\chi_G[S_{n-t}]$ and by (\ref{e5.6.1})
$$\chi_H[S]= \chi_G[S_{n-t}]-A_{0,n-t+1}p^{n-t}+A_{t,n}p^n$$ $$+\sum_{i=1}^{t-1}{(A_{i,n-t+i}-A_{i,n-t+i+1})p^{n-t+i}}+\delta_S(H^\Delta)p^n$$
$$=-a_1-a_{n-t+1}p^{2n-2t-1}+\sum_{i=2}^{n-t}{a_ip^{n-t+i-1}}-\sum_{i=2}^{t}{a_ip^{n-t+i-1}}$$
\begin{equation} \label{e5.10.1}-\sum_{i=2}^{n-t}{a_ip^{2i-2}}-\sum_{i=2}^{n-t}{a_ip^{2i-3}}+\delta_ S(H^\Delta)p^n.
\end{equation}

If $S$ is a $(0,\dots ,0,a_n)$-homogeneous set, then $\chi_F[S]=0$ whenever $l(F)=n-1$ and $\chi_H[S]=-a_np^{n-1}+p^n\delta_S(H^\Delta)$ whenever $l(H)=n$. Therefore $a_n=0$. If $S$ is a $(p+1,p-1,\dots ,p-1,a_n)$-homogeneous set, then we turn to a complement.

In other cases $|\{ \chi_H[S]\mid l(H)=n-1\} |=2$ and $|\{ \chi_H[S]\mid l(H)=n\} |=2$. Therefore $S$ being an SRCG implies $x_{n-1}[S]=x_n[S]$. Since $x_n[S]=-\sum_{i=1}^{n}{a_ip^{i-1}}$, $x_{n-1}[S]=-a_1-a_2p-
\sum_{i=2}^{n-1}{a_ip^i}$, we have 
\begin{equation} \label{e5.10.2}
a_2p+\sum_{i=2}^{n-1}{a_ip^i}=\sum_{i=2}^{n}{a_ip^{i-1}}.
\end{equation}
In the last equality every degree of $p$ from 1 to $n-1$ occurs in the left and in the right hands exactly once with a coefficient $a_i$, $0\leq a_i\leq p-1$, and for all $i>2$ there exists $j<i$ such that $a_i=a_j$. Consequently, $a_2=\cdots =a_n$.

Conversely, if we assume that $a_2=\cdots =a_n$, then substituting into (\ref{e5.10.1}) we obtain $\chi_H[S]=-\sum_{i=1}^{n}{a_ip^{i-1}}+\delta_S(H^\Delta)p^n$ for $H\in \nabla_t\backslash \nabla_{t-1}$. We remark that in this case $k=x_0[S]=s-sp^n$. $\square$

\section{Non-homogeneous strongly regular Cayley \\ graphs over ${\mathbb Z}_{p^n} \oplus {\mathbb Z}_{p^n}$}

If $n=1$, then each subset of standard basis elements of $W(G)$ corresponds to an SRCG \cite{golfand-93}. In what follows we assume that $n\geq 2$.

Let us give examples of SRCGs over $G$ which are not defined by a homogeneous set. Take 
$G={\mathbb Z}_4\oplus {\mathbb Z}_4$. The graph which is defined by a union of $(2,0)$-homogeneous set $S$ and the unique block $B\subseteq (\Delta_2\backslash \Delta_1)\backslash S$ is an $(16,10,6,6)$-SRCG and we denote the set of such graphs by $\Gamma_2$. We denote by $\Gamma_2^c$ the set of their complements.  Similarly, take $G={\mathbb Z}_8\oplus {\mathbb Z}_8$. The graph which is defined by a union of $(3,0,0)$-homogeneous set $S$ and all blocks $B\subseteq (\Delta_3\backslash \Delta_2)\backslash S$ is an $(64,45,32,30)$-SRCG and we denote the set of such graphs by $\Gamma_3$. We denote by $\Gamma_3^c$ the set of their complements. A graph from $\Gamma_3^c$ is a $(64,18,2,6)$-SRCG. Strongly regular graphs with these parameters were enumerated in \cite{haemers-01}. 

A complement of a graph which is defined by an $(a_1,\dots ,a_n)$-homogeneous set is a $(p+1-a_1,p-1-a_2,\dots ,p-1-a_n)$-homogeneous set. A homogeneous set which is block equivalent to a complement of an $(a_1,\dots ,a_n)$-homogeneous set is a $(p+1-a_1,p-1-a_2,\dots ,p-1-a_n)$-homogeneous set or $(0,\dots ,0,p-a_i,p-1-a_{i+1},\dots ,p-1-a_n)$-homogeneous set if $a_1=\cdots a_{i-1}=0$ or $(p+1,p-1,\dots ,p-1,p-2-a_i,p-1-a_{i+1},\dots ,p-1-a_n)$-homogeneous set if $a_1=p+1$, $a_2=\cdots =a_{i-1}=p-1$.

\begin{proposition} \label{p6.1}
If  $S$ and $S^h$ define non-trivial SRCGs, then $S=S^h$ or $\Gamma_G(S)\in \Gamma_2$ or  $\Gamma_2^c$ or $\Gamma_3$ or  $\Gamma_3^c$.
\end{proposition}
{\bf Proof. } By Proposition \ref{p5.10} and Corollary \ref{c5.8} if $S^h$ is an $(a_1,\dots ,a_n)$-homogeneous set which defines an SRCG and $S\neq S^h$, then either $a_2=\cdots =a_n =0$ or $a_2=\cdots =a_n =p-1$ . Suppose that $a_2=\cdots =a_n =0$. Then $\chi_G[S]>\chi_G[S^h]$ hence $\chi_G[S]=k''=s+sp^n+p^n+p^{2n}$ and $k''>k'$. Let $B=\bigcup \{ D \mbox{ is a block in }\Delta \mid D\cap S^h=\emptyset\}$. Then $$\chi_G[B]=\sum_{i=2}^{n}{(p^i+p^{i-1}-A_{0,i-1}p)(p^i-p^{i-1})=p^{2n}-p^2+a_1p^2-a_1p^{n+1}}.$$  
Therefore $\chi_G[S]-\chi_G[S^h]\leq p^{2n}-p^2+a_1p^2-a_1p^{n+1}$, where $\chi_G[S^h]=s-sp^n$ and $s=-a_1$, which is equivalent to $a_1-1\geq ((p-2)a_1+1)p^{n-2}$. Hence $n=2$, $p=2$, $a_1=2$ or $n=3$, $p=2$, $a_1=3$. In the first case $S^h\cup B$ defines a graph from $\Gamma_2$ and in the second case $S^h\cup B$ defines a graph from $\Gamma_3$. $\square$

\begin{proposition} \label{p6.2}
Let $H\in \nabla_1\backslash \nabla_0$.  Denote $\Omega_H=(\Delta\backslash \Delta_{n-1})\backslash \Des_{n-1}(H^\Delta)$. Then for each subset $S\subseteq \Delta$ it is fulfilled that $\chi_G[S]-\chi_H[S]=p^n|\Omega_H\cap S|$.
\end{proposition}
{\bf Proof. } The equation is a straightforward consequence of Corollary \ref{c4.4}. $\square$

\begin{proposition} \label{p6.3}
Let $\Gamma_G(S)$ be a non-trivial $(p^{2n},k,\lambda ,\mu )$-SRCG over $G$ with $k=s+sp^n+p^n+p^{2n}$. If $p>2$, then $S$ is a $((p+1)/2,(p-1)/2,\dots ,(p-1)/2)$-homogeneous set which defines an SRCG with Paley parameters. Moreover, $((p+1)/2,(p-1)/2,\dots ,(p-1)/2)$-homogeneous sets exhaust the set of SRCGs with Paley parameters over $G$. If  $p=2$, then $S$ or its complement satisfy $S\backslash S_{n-1}= (S^h\backslash S_{n-1}^h)\cup (\bigcup \{ D \mbox{ is a block in } \Delta\backslash \Delta_{n-1} \mid D\cap S^h=\emptyset\} )$.
\end{proposition}
{\bf Proof. } Let $H\in \nabla_1\backslash \nabla_0$. Then  by Proposition \ref{p6.2} $\chi_G[S]-\chi_H[S]=p^n|\Omega_H\cap S|$. Since $S$ defines a non-trivial SRCG, $\chi_H[S]\in \{ s,s+p^n\}$. By assumption $k=s+sp^n+p^n+p^{2n}$. Therefore 
\begin{displaymath}
|\Omega_H\cap S|=
\left\{ 
\begin{array}{ll}
s+p^n+1 & \mbox{if} \; \chi_H[S]=s, \\
s+p^n & \mbox{if} \; \chi_H[S]=s+p^n. 
\end{array}
\right.
\end{displaymath}
In other words, $|\Omega_H\cap S|=s+p^n+1-\delta_{\{ s+p^n\} }(\chi_H[S])$. Let $S^h$ be an $(a_1,\dots ,a_n)$-homogeneous set. Let $S^*=S_{n-1}\cup (S^h\backslash S_{n-1}^h)$. Then $\Omega_H\cap S^*=\Omega_H\cap S^h$ and $|\Omega_H\cap S^*|=A_{0,n}-(A_{1,n}+\delta_{S^h}(H^\Delta))=-s-\delta_{S^h}(H^\Delta)$.

If $0<a_n<p-1$, then $S^*=S$, $\delta_{S^h}(H^\Delta)=\delta_{\{ s+p^n\} }(\chi_H[S])$ whence $s+p^n+1=-s$. Therefore $p$ is an odd prime and $$-s=(p^n+1)/2\equiv (p+1)/2(\mod \; p) \not\equiv 1(\mod \; p).$$ Thus the first coordinate of $S^h$ is $(p+1)/2$ and $s$ defines the homogeneous coordinates of $S^h$ which are equal to $((p+1)/2,(p-1)/2,\dots ,(p-1)/2)$. By Corollary \ref{c5.8} $S=S^h$ since $(p-1)/2\neq 0$ and  $(p-1)/2\neq p-1$. Thus $S$ is homogeneous and defines an SRCG with Paley parameters by Proposition \ref{p5.10}.

Assume now that $S\subseteq \Delta$ is a set which defines an SRCG with Paley parameters. Then $\lambda -\mu =r+s=s+p^n+s=-1$. Using arguments from the previous paragraph we obtain that SRCGs with Paley parameters over $G$ are $((p+1)/2,(p-1)/2,\dots ,(p-1)/2)$-homogeneous sets.

Consider the case of $a_n=0$ and $p>2$. Let $B=\bigcup \{ D \mbox{ is a block in }\Delta \mid  D\subseteq \Omega_H, D\cap S^h=\emptyset \} $ and $H\in \nabla_1\backslash \nabla_0$. Then $|B|=p^n-(A_{1,n-1}p+a_1-\delta_{S^h}(H^\Delta))p$ and
$$0=|\Omega_H\cap S|-|\Omega_H\cap S^*|-|B\cap S|\geq |\Omega_H\cap S|-|\Omega_H\cap S^*|-|B|$$
$$=s+p^n+1-\delta_{\{ s+p^n\} }(\chi_H[S])+s+\delta_{S^h}(H^\Delta)-p(p^{n-1}-A_{1,n-1}p)+a_1p-\delta_{S^h}(H^\Delta)p$$
\begin{equation} \label{e6.3.1}
=\sum_{i=1}^{n-1}{a_i(p^i-2p^{i-1})}+1-\delta_{\{ s+p^n\} }(\chi_H[S])-(p-1)\delta_{S^h}(H^\Delta). 
\end{equation}
If $p>2$ and $a_i\neq 0$ for some $2\leq i<n$, then (\ref{e6.3.1}) is strongly positive. Therefore $a_2=\cdots =a_{n-1}=a_n=0$. Now by Proposition \ref{p6.1} $S=S^h$ and $k=s-sp^n$. Together with $k=s+sp^n+p^n+p^{2n}$ this implies $s=-(p^n+1)/2$ and contradicts to $a_n=0$.

Consider now the case of $a_n=0$ and $p=2$. Let $H\in \nabla_1\backslash \nabla_0$ and 
$B=\bigcup \{ D \mbox{ is a block in }\Delta \mid D\subseteq \Omega_H,D\cap S^h=\emptyset\} $. Then $$0=|\Omega_H\cap S|-|\Omega_H\cap S^*|-|B\cap S|\geq |\Omega_H\cap S|-|\Omega_H\cap S^*|-|B|$$
$$=1-\delta_{\{ s+p^n\} }(\chi_H[S])-\delta_{S^h}(H^\Delta).$$
Since $D\cap S^h=\emptyset $ for each block $D\subseteq B$, it holds that $|B\cap S|$ and $|B|$ are divisible by $p$. Therefore 
$\delta_{\{ s+p^n\} }(\chi_H[S])+\delta_{S^h}(H^\Delta)=1$ and $B\cap S=B$. Since $\Delta\backslash \Delta_{n-1}=\bigcup_{H\in \nabla_1\backslash \nabla_0}{\Omega_H}$, it holds that $\bigcup \{ D \mbox{ is a block in } \Delta\backslash \Delta_{n-1} \mid D\cap S^h=\emptyset\}\subseteq S$. In particular, graphs from $\Gamma_2$ and $\Gamma_3$ satisfy this condition.

If $\Gamma _G(S)$ is an SRCG with the valency $k=s+sp^n+p^n+p^{2n}$, then its complement has the valency of the same type. Then the case of $a_n=p-1$ is complement to the case of $a_n=0$. $\square$

\begin{proposition} \label{p6.4}
Let $S$ define a non-trivial SRCG over $G$. If $k=s-sp^n$, then $S\backslash S_{n-1}= S^h\backslash S_{n-1}^h$.  
\end{proposition}
{\bf Proof. } Let $H\in \nabla_1\backslash \nabla_0$. Let $S^*=S_{n-1}\cup (S^h\backslash S_{n-1}^h)$. Then 
$$\chi_G[S^*]=\chi_G[S_{n-1}]+A_{0,n}(p^n-p^{n-1}),$$
$$\chi_H[S^*]=\chi_G[S_{n-1}]+(A_{1,n}+\delta_{H^\Delta}(S^h))p^n-A_{0,n}p^{n-1}.$$
Therefore by Proposition \ref{p6.2} $$|\Omega_H\cap S^*|=A_{0,n}-A_{1,n}-\delta_{H^\Delta}(S^h)=-s-\delta_{H^\Delta}(S^h).$$ Again by Proposition \ref{p6.2} and by $k=s-sp^n$ we have $|\Omega_H\cap S|=-s-\delta_{\{ s+p^n\} }(\chi_H[S])$. By construction of $S^h$ either $S^*\subseteq S$ or $S\subseteq S^*$. Therefore 
either $\Omega_H\cap S^*\subseteq \Omega_H\cap S$ or $\Omega_H\cap S\subseteq \Omega_H\cap S^*$. From which it follows that $|(\Omega_H\cap S)\ominus (\Omega_H\cap S^*)|=||\Omega_H\cap S|-|\Omega_H\cap S^*||\in \{ 0,1\}$, where $\ominus $ denotes the symmetric difference. Since $\Omega_H\cap S$ and $\Omega_H\cap S^*$ are block equivalent, the cardinality of their symmetric difference is divisible by $p$. Therefore $|\Omega_H\cap S|=|\Omega_H\cap S^*|$ and $\Omega_H\cap S=\Omega_H\cap S^*$. Since $\Delta\backslash \Delta_{n-1}=\bigcup_{H\in \nabla_1\backslash \nabla_0}{\Omega_H}$, we obtain $(\Delta\backslash \Delta_{n-1})\cap S=(\Delta\backslash \Delta_{n-1})\cap S^*$ from which it follows that $S\backslash S_{n-1}=S^h\backslash S_{n-1}^h$. $\square$ 

\begin{proposition} \label{p6.5}
Let $S$ define a non-trivial SRCG over $G$, $p>2$. Then either $S$ is an $(a_1,a_2,\dots ,a_2)$-homogeneous set or $S^h$ is an $(a_1,\dots ,a_n)$-homogeneous set with $a_{n-1}=0$ or $a_{n-1}=p-1$.
\end{proposition}
{\bf Proof. } Now we can assume that $n\geq 3$. According to Corollary \ref{c5.9} there exists a block set $U$ such that $S=S^h\ominus U$. Denote $R=S^h\cap U$, $T=U\backslash S^h$. Let $F\in \nabla_{n-1}\backslash \nabla_{n-2}$. Then $\chi_F[S]=\chi_F[S^h]+\chi_F[T]-\chi_F[R]$. We have by Proposition \ref{p5.6}  $|\{ \mbox{ residue of  }\chi_F[S^h] \mbox{ modulo }p^n\mid l(F)=n-1\} |=1$. Furthermore, 
$|\{ \mbox{ residue of  }\chi_F[S] \mbox{ modulo }p^n \mid l(F)=n-1\} |=1$.
Therefore 
$|\{ \mbox{ residue of  }(\chi_F[T]-\chi_F[R] )\mbox{ modulo }p^n\mid l(F)=n-1\} |=1$. By Proposition \ref{p6.4} $S\backslash S_{n-1}=S^h\backslash S_{n-1}^h$. Then $T^1\subseteq \Delta _{n-2}^1$, $R^1\subseteq \Delta _{n-2}^1$.
Denote ${\bar R}=\bigcup_{\{ m \mid a_m=p-1,2\leq m\leq n-1\} }{(\Delta_m\backslash \Delta_{m-1})}\backslash R$. Since $\chi_F[\Delta_m\backslash \Delta_{m-1}]=-p^2\delta_{\{ 2\} }(m)$ for $n\geq m \geq 2$, it follows that $\chi_F[{\bar R}]=-p^2\delta_{\{ m\mid a_m=p-1\} }(2)-\chi_F[R]$. Denote $Q=T\cup {\bar R}$ and $\rho=\delta_{\{ m \mid a_m=p-1\} }(2)$. By Proposition \ref{p5.7} $Q^1$ is block equivalent to the unique $(b_1,b_2,\dots ,b_{n-2})$-homogeneous set $Q^{1h}\subseteq \Delta_{n-2}^1$ with $Q^{1h}\cap \Delta_1^1=Q^1\cap \Delta_1^1$. Therefore
$$\chi_F[Q]=\chi_F[Q^{1h}]=x_{n-1}[Q]+\delta_{Q^{1h}}(\Sons(\Father(F^\Delta)))p^n.$$ 
Thus if $Q^{1h}$ is not a $(p+1,p-1,\dots ,p-1)$-homogeneous set, then
$$\chi_F[S]=s+\delta_{\{ s+p^n\} }(\chi_F[S])p^n=\chi_F[S^h]+\chi_F[T]-\chi_F[R]
=\chi_F[S^h]+\chi_F[Q]+\rho p^2$$ 
\begin{equation}\label{e6.5.2}
=x_{n-1}[S^h]+\delta_{S^h}(F^\Delta)p^n+x_{n-1}[Q]+\delta_{Q^{1h}}(\Sons(\Father(F^\Delta)))p^n+\rho p^2. 
\end{equation}
If $Q^{1h}$ is a $(p+1,p-1,\dots ,p-1)$-homogeneous set, then $\chi_F[S]=s+\delta_{\{ s+p^n\} }(\chi_F[S])p^n$ $=x_{n-1}[S^h]+\delta_{S^h}(F^\Delta)p^n-p^2+\rho p^2 $.

Consider the cases when $Q^{1h}$ is a (0,\dots ,0)-homogeneous set and $\rho =0$ or $Q^{1h}$ is a $(p+1,p-1,\dots ,p-1)$-homogeneous set and $\rho =1$. 
If $a_{n-1}=0$ or $a_{n-1}=p-1$ then there is nothing to prove.
If $a_{n-1}\neq 0$ and $a_{n-1}\neq p-1$ then $\delta_{S^h}(F^\Delta)$ has two values and $s=x_{n-1}[S^h]$. Now, $S^h$ is an $(a_1,a_2,\dots ,a_2)$-homogeneous set by (\ref{e5.10.2}). Then by Proposition \ref{p6.1} $S=S^h$.

If $Q^{1h}$ is a $(p+1,p-1,\dots ,p-1)$-homogeneous set and $\rho =0$, then $a_2\neq p-1$, $Q\cap (\Delta_2\backslash \Delta_1)=T\cap (\Delta_2\backslash \Delta_1)$ and by Corollary \ref{c5.8}  $a_2=0$, $a_1=0$ and $s=x_{n-1}[S^h]-p^2$. Thus $\sum _{i=3}^{n-1}a_ip^{i-1}=(a_np^{n-1}-p^2)/(p-1)$ and $$\sum _{i=3}^{n-1}a_ip^{i-1}=\sum _{i=3}^{n-1}a_np^{i-1}+p^2(a_n-1)/(p-1).$$ This implies that $a_n=1$ and $\sum _{i=3}^{n-1}a_ip^{i-2}=\sum _{i=3}^{n-1}p^{i-2}$.
Therefore $S^h$ is an $(0,0,1,\dots ,1)$-homogeneous set and $S=S^h\cup (\Delta _2\backslash \Delta _1)$. The equality $k=s-sp^n=\chi _G[S^h]+\chi _G[\Delta _2\backslash \Delta _1]$ implies $n=4$. Thus $S^h$ is a $(0,0,1,1)$-homogeneous set. By direct calculations one can check that $s\neq x_2[S]$.

If $Q^{1h}$ is a $(0,\dots ,0)$-homogeneous set and $\rho =1$, then $a_2=p-1$, $a_1=p+1$, and $s=x_{n-1}[S^h]+p^2$. Analogously to the previous case, $S^h$ is a $(p+1,p-1,p-2,\dots ,p-2)$-homogeneous set. Turning to a complement we obtain that there is no non-trivial SRCG over $G$ which is defined by a set $S$ such that $S^h$ is  $(p+1,p-1,p-2,\dots ,p-2)$-homogeneous.

The fact that $Q^1$ is block equivalent to a homogeneous set which is neither $(0,\dots ,0)$-homogeneous nor $(p+1,p-1,\dots ,p-1)$-homogeneous implies by Corollary \ref{c5.8} that $a_{n-1}=0$ or $a_{n-1}=p-1$. $\square$

\begin{lemma} \label{l6.6}
Let $S$ define a non-trivial SRCG over $G$, $p>2$. Then either $S$ is an $(a_1,a_2,\dots ,a_2)$-homogeneous set or one of the sets $S^h$ or $(\Delta \backslash \Delta_0)\backslash S^h$ is an $(a_1,0,\dots ,0,a_n)$-homogeneous set.
\end{lemma}
{\bf Proof. } We shall use the notations of the previous proposition. According to Proposition \ref{p6.5} we can assume that $n\geq 4$ and since the case of $a_{n-1}=p-1$ is complement to $a_{n-1}=0$, we assume $a_{n-1}=0$. This assumption entails that if $l(F)=n-1$ and $\delta_{S^h}(F^\Delta)=1$, then $\delta_{Q^{1h}}(\Sons(\Father(F^\Delta)))=0$ because a block cannot be included in $S$ partially and completely at the same time. We shall use this consideration in the proof.

Assume the contrary. Then there exists $j$, $2\leq j\leq n-2$, such that $a_j\neq 0$ and $a_i=0$ for each $i$, such that $j<i\leq n-1$. Then $\delta_{S^h}(F^\Delta)$ has two values  on $\nabla_{n-1}\backslash \nabla_{n-2}$ and, consequently, (\ref{e6.5.2}) implies $s=x_{n-1}[S^h]+x_{n-1}[Q]+\rho p^2$. We have $s=-\sum_{i=1}^{n}{a_ip^{i-1}}$, $ x_{n-1}[S^h]=-a_1-a_2p-\sum_{i=2}^{n-1}{a_ip^i}$ by (\ref{e5.10.1}) which implies that 
\begin{equation}\label{e6.6.1}
0\leq -x_{n-1}[Q]=\sum_{i=1}^{n-2}{b_ip^{i+1}}=a_np^{n-1}-(p-1)\sum_{i=3}^{n-1}{a_ip^{i-1}}-a_2p^2+\rho p^2.
\end{equation}

Then (\ref{e6.6.1}) implies $a_n>0$ and
\begin{equation} \label{e6.6.2}
b_{n-2}\leq a_n-1.
\end{equation}

It holds that $$p^{n-1}+p^2\geq \sum_{i=1}^{n-3}{b_ip^{i+1}}
=(a_n-b_{n-2})p^{n-1}-(p-1)\sum_{i=3}^{j}{a_ip^{i-1}}-a_2p^2+\rho p^2$$
$$\geq (a_n-b_{n-2})p^{n-1}-p^{j+1}+p^j+p^2.$$

These inequalities and (\ref{e6.6.2}) imply
\begin{equation} \label{e6.6.4}
b_{n-2}=a_n-1.
\end{equation}

Further, $(1-p)\sum_{i=3}^{j+1}{a_ip^{i-1}}\equiv a_2p^2-\rho p^2+\sum_{i=1}^{j-1}{b_ip^{i+1}}(\mod \; p^{j+1})$.

If $b_1=b_2=\cdots =b_{j-1}=0$ and $\rho=0$ or $b_1=p+1$, $b_2=\cdots =b_{j-1}=p-1$ and $\rho=1$,
then  $(1-p)\sum_{i=3}^{j+1}{a_ip^{i-1}}\equiv a_2p^2(\mod \; p^{j+1})$, $\sum_{i=3}^{j+1}{a_ip^{i-1}}\equiv  a_2p^2+\dots +a_2p^j(\mod \; p^{j+1})$.
Therefore $a_2=a_3=\cdots =a_j=a_{j+1}=0$ which contradicts to $a_j\neq 0$.

If $b_1=p+1$, $b_2=\cdots =b_{j-1}=p-1$ and $\rho=0$, then $a_2=0$ and 
$(1-p)\sum_{i=3}^{j+1}{a_ip^{i-1}}\equiv p^2(\mod \; p^{j+1})$,  $a_{j+1}=1$ which contradicts to $a_{j+1}=0$. 

If $b_1=b_2=\cdots =b_{j-1}=0$ and $\rho=1$, then $a_2=p-1$ and  $(1-p)\sum_{i=3}^{j+1}{a_ip^{i-1}}\equiv (p-2)p^2(\mod \; p^{j+1})$,  $a_{j+1}=p-2$ which contradicts to $a_{j+1}=0$. 

Therefore $b_1=b_2=\cdots =b_{j-1}=0$ contradict to the assumption and $b_1=p+1$, $b_2=\cdots =b_{j-1}=p-1$ contradict to the assumption. By Corollary \ref{c5.8} 
\begin{equation} \label{e6.6.3}
a_j=p-1.
\end{equation}

Thus 
\begin{equation} \label{e6.6.31}
\sum_{i=1}^{n-3}{b_ip^{i+1}}=p^{n-1}-(p-1)\sum_{i=3}^{j}{a_ip^{i-1}}-a_2p^2+\rho p^2\geq 
p^{n-1}- p^{j+1}+ p^j+p^2.
\end{equation}
If $j<n-2$, then the last inequality gives $b_j=p-1$; moreover, $b_1=p+1$, $b_2=\cdots =b_{j-1}=p-1$ contradict to the assumption. Thus $b_j=p-1$ and $a_{j+1}=0$ imply that there exists a block $B\subseteq \Delta _j^1\backslash \Delta_{j-1}^1$ such that $p-1$ of its vertices are included in $T^1$. On the other hand, (\ref{e6.6.3}) and $a_{j+1}=0$ imply that each block in $\Delta _j^1\backslash \Delta_{j-1}^1$ includes at least $p-1$ vertices which have only one vertex from $S$. Since $p>2$, we have a contradiction.

If $j=n-2$, then $b_1=b_2=\cdots =b_{n-3}=0$ contradict to the assumption and $b_1=p+1$, $b_2=\cdots =b_{n-3}=p-1$ contradict to the assumption. If, in addition, $a_n>1$, then by (\ref{e6.6.4}) $b_{n-2}\geq 1$ and there exists a block $B\subseteq \Delta_{n-2}^1 \backslash \Delta_{n-3}^1$ such that $1<|B\cap T^1|<p$. Then we have a contradiction analogously to the previous case.

Now consider the case of $j=n-2$ and $a_n=1$. Then $a_{n-2}=p-1$ by (\ref{e6.6.3}), $a_{n-1}=0$, $b_{n-2}=0$. If $a_{n-2}=p-1$ and $|B\cap S^h|=p-1$ for a block $B\subseteq \Delta _{n-2}\backslash \Delta_{n-3}$, then $B^1\in Q^1$. If $b_{n-3}>0$, then $B^1\in Q^{1h}$. If $b_{n-3}=0$, then $b_1=p+1$, $b_2=\cdots =b_{n-4}=p-1$ by (\ref{e6.6.31}) and also $B^1\in Q^{1h}$. Therefore $|\Sons(B^1)\cap Q^{1h}|=b_{n-2}+1=1$. On the other hand, if $a_{n-1}=0$ and $|B\cap S^h|=p-1$, then there exists a unique block $C$ such that $C^1\in \Sons(B^1)$ and $C\cap S^h=\emptyset$. 
Then $\Sons(B^1)\cap Q^{1h}=\{ C^1\in \Sons(B^1) \mid \Father(C)\not\in S^h\} $. If $a_{n-2}=p-1$ and $|B\cap S^h|=p$ for a block $B\subseteq \Delta _{n-2}\backslash \Delta_{n-3}$, then $B\subseteq S^{h}$ and since $a_{n-1}=0$, there is no block $C^1\in \Sons(B^1)$ such that $C\cap S^h=\emptyset$.
Then $\Sons(B^1)\cap Q^1=\emptyset$ and $\Sons(B^1)\cap Q^{1h}=\{ C^1\in \Sons(B^1) \mid \Father(C)\not\in S^h\} =\emptyset$ since $b_{n-2}=0$.  Therefore $Q^{1h}\cap \Delta _{n-2}^1\backslash \Delta _{n-3}^1=\{ C^1 \in \Delta _{n-2}^1\backslash \Delta _{n-3}^1 \mid \Father(C)\not\in S^h\} $.
This implies that  $Q^{1h}$ has coordinates of a homogeneous set which is block equivalent to a complement of $S_{n-2}^h$. Following the paragraph before Proposition \ref{p6.1} we set 
\begin{displaymath}
t=
\left\{ 
\begin{array}{ll}
0 & \mbox{if} \; Q^{1h}  \mbox{ is a } (p+1-a_1,p-1-a_2,\dots ,p-1-a_{n-2})\\
 & \mbox{-homogeneous set,}\\
-1 & \mbox{if} \; Q^{1h}  \mbox{ is a } (0,\dots ,0,p-a_i,p-1-a_{i+1},\dots ,p-1-a_{n-2})\\
 & \mbox{-homogeneous set,}\\
1 & \mbox{if} \; Q^{1h}  \mbox{ is a } (p+1,p-1,\dots ,p-1,p-2-a_i,p-1-a_{i+1},\\
 & \dots ,p-1-a_{n-2}) \mbox{-homogeneous set.}
\end{array}
\right.
\end{displaymath}
By definition $t=-1$ implies $a_1=0$ and $t=1$ implies $a_1=p+1$.
Since $a_{n-1}=0$ and $a_n=1$, $$-x_{n-1}[Q]=p^{n-1}-(p-1)\sum_{i=3}^{n-2}{a_ip^{i-1}}-a_2p^2+\rho p^2$$
\begin{equation} \label{e6.6.5} =(p+1-a_1)p^2+\sum_{i=2}^{n-2}{(p-1-a_i)p^{i+1}}+tp^2.
\end{equation}
Then $(p^2-p+1)\sum_{i=3}^{n-2}{a_ip^{i-1}}=p^n-p^{n-1}-a_1p^2-a_2p^3+a_2p^2+p^2-\rho p^2+tp^2$.
Since $a_{n-2}=p-1$, in the case of $n \geq 5$ we can rewrite the last equality as follows: 
$$(p^2-p+1)\sum_{i=3}^{n-3}{a_ip^{i-1}}=$$
\begin{equation} \label{e6.6.6}
p^{n-1}-2p^{n-2}+p^{n-3}-a_1p^2-a_2p^3+a_2p^2+p^2-\rho p^2+tp^2.
\end{equation}
If $a_1=0$, then $t\neq 1$. If $a_1=p+1$, then $t\neq -1$. This implies
$$(-p^2+p-2)p^2\leq -a_1p^2-a_2p^3+a_2p^2+p^2-\rho p^2+tp^2\leq p^2.$$
Then 
$$p^{n-3}-p^{n-4}-p^{n-5}-p^2+\frac{p^{n-5}-p^2}{p^2-p+1}$$
$$\leq \sum_{i=3}^{n-3}{a_ip^{i-1}}\leq p^{n-3}-p^{n-4}-\frac{p^{n-3}-p^{n-4}-p^2}{p^2-p+1}.$$
Therefore if $n\geq 6$, then $a_{n-3}\leq p-2$. If $n\geq 7$, then $a_{n-3}=p-2$. If $n=6$, then $(p^2-p+1)(p-a_3)=p^2+a_1+a_2p-a_2-1+\rho -t\leq 2p^2-p+3$ by (\ref{e6.6.6}) and $a_{n-3}=p-2$. Hence if $p>2$, then either $b_1=b_2=\cdots =b_{n-4}=0$ or 
$b_1=p+1$, $b_2=\cdots =b_{n-4}=p-1$ by Corollary \ref{c5.8}. This implies that  $a_1=a_2=\cdots =a_{n-4}=0$ or $a_1=p+1$, $a_2=\cdots =a_{n-4}=p-1$. Substituting these values into  (\ref{e6.6.6}) we obtain a contradiction.

If $n=5$, then $S^h$ is an $(a_1,a_2,p-1,0,1)$-homogeneous set and (\ref{e6.6.6}) implies $a_1+a_2(p-1)-1+\rho -t=(p-1)^2$. Therefore if $p>3$, then $p-2\leq a_2\leq p-1$.

Consider the case of $n=5$, $p=3$, $a_2=0$. Then $\rho =0$, $a_1=5+t$ which imply $t=-1$ and $a_1=p+1$, but $t=-1$ implies $a_1=0$, a contradiction.

Consider the case of $n=5$, $a_2=p-2$. Then $a_1=p+t$. Therefore $t=1$, $a_1=p+1$, $b_1=p+1$  or $t=0$, $a_1=p$, $b_1=1$ which contradict to $a_2=p-2$ by Corollary \ref{c5.8}.

Consider the case of $n=5$, $a_2=p-1$. In this case $\rho =1$, $0\leq a_1=t$. If $t=1$, then $a_1=p+1$, but $a_1=t=1$ and we have a contradiction. If $t=a_1=0$, then $S^h$ is a $(0,p-1,p-1,0,1)$-homogeneous set, $Q^{1h}$ is a $(p+1,0,0)$-homogeneous set which is complement to $S_{n-2}^h$.

Consider the case of $n=4$. Then $S^h$ is an $(a_1,p-1,0,1)$-homogeneous set and $\rho =1$. By 
(\ref{e6.6.5}) we have $a_1=p+t-1$. If $t=-1$ or $t=1$, we have a contradiction to the definition of $t$. If $t=0$, then $S^h$ is a $(p-1,p-1,0,1)$-homogeneous set, $Q^{1h}$ is a (2,0)-homogeneous set. 

A construction of $S$ in the last 2 cases leads to a single set in each case (up to automorphism of the tree $\Delta$) and straightforward computation of the principal character leads to a contradiction with $k=s-sp^n$.
 $\square$

{\bf Proof of Theorem \ref{t1.6}. } If $S$ is not an $(a_1,a_2,\dots ,a_2)$-homogeneous set, then by Lemma \ref{l6.6} $S^h$ is an $(a_1,0,\dots ,0,a_n)$-homogeneous set with $a_n>0$ or its complement. Let $S^h$ be an $(a_1,0,\dots ,0,a_n)$-homogeneous set. 
Then  $\chi_H[S^h]=-a_1+\delta _{S^h}(H^\Delta )p^n$ for $H\in \nabla _{n-1}\backslash \nabla _1$ by (\ref{e5.10.1}). If a block $B\subseteq \nabla _{n-1}\backslash \nabla _1$ satisfies $|\{ \chi_H[S^h]| H\in B\}|=1$, then $\chi_H[S^h]=-a_1$ for $H\in B$. $S\backslash S^h$ is a block set in $\Delta _{n-1}\backslash \Delta _1$ by Proposition \ref{p6.4} since $S$ has a principal eigenvalue $k=s-sp^n$ by Proposition \ref{p6.3}. For all block $B\subseteq \nabla _{n-1}\backslash \nabla _1$ it holds that $|\{ \chi_H[S\backslash S^h]| H\in B\}|=1$ since $S\backslash S^h$ is a block set. Therefore by Theorem \ref{t2.1}  $Q=\varphi (S\backslash S^h)$ defines an SRCG over $\varphi (pG)$ with non-principal eigenvalues $-a_np^{n-3}$ or $-a_np^{n-3}+p^{n-2}$. If $n>3$, then $0$ and $-1$ are not eigenvalues of $Q$. Hence $Q$ is a non-trivial SRCG. Thus $Q^h$ is a $(0,\dots ,0,a_n)$-homogeneous set or a $(p,p-1,\dots ,p-1,a_n-1)$-homogeneous set. If $n=3$ then $Q$ is an $(a_3)$-homogeneous set. If $a_3>1$, then an $(a_3)$-homogeneous set defines a non-trivial SRCG over ${\mathbb Z}_p \oplus {\mathbb Z}_p$.

According to Proposition \ref{p6.3} $x_0[Q]=-a_np^{n-3}(1-p^{n-2})$. In addition, $x_1[S\backslash S^h]=x_0[S\backslash S^h]=x_0[Q]p^2$ and $x_1[S^h]=-a_1-a_np^{2n-3}$. Then $x_1[S]=x_1[S^h]+x_1[S\backslash S^h]$ which implies that the graphs mentioned in the theorem are strongly regular. $\square$

\section{Corollaries}

\begin{corollary} \label{c7.1}
A non-trivial SRCG over ${\mathbb Z}_{p^2} \oplus {\mathbb Z}_{p^2}$ is either defined by an $(a_1,a_2)$-homogeneous set, where $(a_1,a_2)\not\in \{(1,0),(p,p-1),(0,0),(p+1,p-1)\}$, or it is the Clebsch graph from $\Gamma_2^c$ or its complement from $\Gamma_2$.
\end{corollary}

\begin{theorem} [\cite{leung-95}] \label{t7.2}
Suppose that there exists an SRCG with Paley parameters over a finite abelian group $A$ of rank 2. Then $A$ is isomorphic to ${\mathbb Z}_{p^n} \oplus {\mathbb Z}_{p^n}$, where $p$ is an odd prime and $n$ is a positive integer.
\end{theorem}

\begin{corollary} \label{c7.3}
SRCGs with Paley parameters $(\nu ,(\nu -1)/2,(\nu -5)/4,(\nu -1)/4)$ over a finite abelian group of rank 2 are defined by $((p+1)/2,(p-1)/2,\dots ,(p-1)/2)$-homogeneous sets over  ${\mathbb Z}_{p^n} \oplus {\mathbb Z}_{p^n}$, where $p$ is an odd prime.
\end{corollary}

\begin{remark} \label{r7.4}
Let $\{e\}\cup A_1\cup \dots \cup A_d$ be a partition of the group ${\mathbb Z}_{p^n} \oplus {\mathbb Z}_{p^n}$ such that each $A_i$, $1\leq i\leq d$, is a set of generators of elements of a homogeneous set which defines an SRCG. Then $\langle 1,\underline{A_1},\dots ,\underline{A_d}  \rangle$ is an S-ring since a disjunctive union of homogeneous sets which define SRCGs is a homogeneous set which defines an SRCG.
\end{remark}

{\bf Acknowledgements}

We would like to thank Ron Adin for helpful remarks. This research was supported by German-Israeli Foundation for Scientific Research and Development.

The full version of this paper was submitted to Discrete Mathematics as "Strongly regular Cayley graphs over the group ${\mathbb Z}_{p^n} \oplus {\mathbb Z}_{p^n}$".


\begin{thebibliography}{99}
\bibitem {bose-63}
R.C. Bose, Strongly regular graphs, partial geometries, and partially balanced designs, Pacific J. Math. 13 (1963) 389-419.
\bibitem {bridges-79}
W.G. Bridges, R.A. Mena, Rational circulants with rational spectra and cyclic strongly regular graphs, Ars Combin. 8 (1979) 143-161.
\bibitem {bridges-82}
W.G. Bridges, R.A. Mena, Rational G-matrices with rational eigenvalues, J. Combin. Theory Ser.A 32(1982) 264-280.
\bibitem{brouwer-89}
A.E. Brouwer, A.M. Cohen, A. Neumaier, Distance-regular  graphs, Springer, Berlin, 1989.
\bibitem{davis-94}
J.A. Davis, Partial difference sets in p-groups, Arch. Math. 63 (1994) 103-110.
\bibitem{delsarte-73}
P. Delsarte, An algebraic approach  to the association schemes of coding theory, Philips Res. Reports (Suppl. 10) (1973)1-97.
\bibitem{faradzev-90}
I.A. Faradzev, A.A. Ivanov, M.H. Klin, Galois correspondence between permutation groups and cellular rings (association schemes), Graphs Combin. 6 (1992) 202-224.
\bibitem{golfand-93}
J.J. Golfand, A.A. Ivanov, M.H. Klin, Amorphic cellular rings,in: I.A. Faradzev, A.A. Ivanov, M.H. Klin, A.J. Woldar (Eds.), Investigations in Algebraic Theory of Combinatorial Objects, Mathematics and Its Applications (Soviet Series), vol. 84, Kluwer Academic Pablishers, Dordrecht 1994,pp. 167-187.
\bibitem{haemers-01}
W.H. Haemers, E. Spence, The pseudo-geometric graphs for generalized quadrangle of order (3,t), European J. Combin. 22(6) (2001) 839-845.
\bibitem{higman-70}
D.G. Higman, A survey of some questions and results about rank 3 permutation groups, Actes Congress Intern. Math. Nice 1970, Vol I (Gauthier-Villars, 1971) 361-365.
\bibitem{kochendorfer-37}
R. Kochendorfer, Untersuchungen {\"u}ber eine Vermutung von W. Burnside, Schr. Math. Sem. Inst. Angew. Math. Univ. Berlin  3 (1937) 155-180.
\bibitem{leung-95}
K.H. Leung, S.L. Ma, Partial difference sets with Paley parameters, Bull. London Math. Soc. 27 (1995) 553-564.
\bibitem{ma-89}
S.L. Ma, On association schemes, Schur rings, strongly regular graphs and partial difference sets, Ars Combin. 27 (1989) 211-220.
\bibitem{muzychuk-87}
M.E. Muzychuk, V-rings of permutation groups with invariant metric, Ph.D Thesis, Department of Mathematics, University of Kiev, 1987 (in Russian).
\bibitem{seidel-69}
J.J. Seidel, Strongly regular graphs, in: W.T. Tutte, (Ed.), Recent Progress in Combinatorics, Academic Press, New York, 1969, pp. 185-197.
\bibitem{wielandt-64}
H. Wielandt, Finite Permutation Groups, Academic Press, New York, 1964.

\end{thebibliography}
\end{document}